\documentclass[12pt]{amsart}

\usepackage{amsfonts}
\usepackage{amsmath}
\usepackage{amssymb}
\usepackage{latexsym}

\newtheorem{theorem}{Theorem}
\theoremstyle{plain}

\newtheorem{definition}[theorem]{Definition}
\newtheorem{proposition}[theorem]{Proposition}

\newtheorem{remark}[theorem]{Remark}

\numberwithin{equation}{section}
\renewcommand\endproof{\hfill $\Box$\vskip 0.15in}

\def\pr{\noindent {\bf {Proof.\ }}}

\newcommand\bel[1]{\begin{equation}\label{#1}}
\newcommand\ee{\end{equation}}

\newcommand\und[1]{{\tt{#1}}}

\newcommand{\ol}[1]{{\overline{#1}}}
\newcommand{\wt}[1]{{\widetilde{#1}}}

\newcommand{\bD}{{\mathbb{D}}}

\newcommand{\bR}{{\mathbb{R}}}
\newcommand{\bE}{{\mathbb{E}}}
\newcommand{\bP}{{\mathbb{P}}}

\newcommand{\bA}{{\mathbf{V}}}

\newcommand{\bH}{{\mathbf{H}}}

\newcommand{\HX}{{\mathbb{H}}_{\mfX}}
\newcommand{\HXt}{\HX^{(1)}}
\newcommand{\HB}{{\mathbb{H}}_{\mfB}}
\newcommand{\HBt}{{\HB^{(1)}}}

\newcommand{\cF}{{\mathcal{F}}}
\newcommand{\cI}{{\mathcal{I}}}
\newcommand{\cJ}{{\mathcal{J}}}
\newcommand{\cK}{{\mathcal{K}}}

\newcommand{\cR}{{\mathcal{R}}}

\newcommand{\mfX}{{\mathfrak{X}}}
\newcommand{\mfB}{{\mathfrak{B}}}

\def\Ito{It\^{o}}

\begin{document}

\today


\title[Stochastic Integration]
{From Random Processes to Generalized Fields: A Unified Approach to Stochastic
Integration}
\author{S. V. Lototsky}
\curraddr[S. V. Lototsky]{Department of Mathematics, USC\\
Los Angeles, CA 90089}
\email[S. V. Lototsky]{lototsky@math.usc.edu}
\urladdr{http://www-rcf.usc.edu/$\sim$lototsky}
\author{K. Stemmann}
\curraddr[K. Stemmann]{Department of Mathematics, USC\\
Los Angeles, CA 90089}
\email[K. Stemmann]{stemmann@usc.edu}
\thanks{S. V. Lototsky  acknowledges  support from
 the Sloan Research Fellowship, the NSF
CAREER award DMS-0237724, and the ARO Grant DAAD19-02-1-0374.
The work of K. Stemmann was partially supported by the NSF Grant
DMS-0237724}

\begin{abstract}
The paper studies stochastic integration with respect to Gaussian processes and
fields. It is more convenient to work with a field than a process:
 by definition, a field is
  a collection of stochastic integrals for a class of deterministic integrands.
  The problem is then to extend the definition to random integrands.
 An orthogonal decomposition of
chaos space of the random field leads to two such extensions,
corresponding to the \Ito-Skorokhod and the Stratononovich integrals,
  and provides an efficient tool to
study these integrals, both analytically and numerically.
For a Gaussian process, a natural definition of the integral follows
from  a canonical
correspondence between random processes and a special class of random
fields.
\end{abstract}

\maketitle

\section{Introduction}
While stochastic integral with respect to a standard Brownian
motion is a well-studied object,  integration
with respect to other Gaussian processes is currently an area of active research,
and the  fractional Brownian motion is receiving most of the attention
\cite[etc.]{AMN1,A,DH,DU,DHP,KBR,Lin,PT}
The objective of this paper is
to define and investigate stochastic
integrals with respect to arbitrary
Gaussian processes and fields using chaos expansion.
The motivation comes from the paper by Al\`{o}s et al. \cite{AMN}
and the book by P. Major \cite{Maj}.

In  \cite{AMN}, the authors study stochastic integration with respect to
the Gaussian process $\int_0^t K(t,s)dW(s)$, where $K$ is a suitable kernel
function and $W$ is a standard Brownian motion. In  \cite{Maj},
the author studies stochastic integration with respect to generalized Gaussian
fields. While \cite{AMN} and \cite{Maj} pursue different goals and
work with  different objects,  generalized fields and the chaos expansion,
appearing in both \cite{AMN} and \cite{Maj}, are the
unifying ideas.

A generalized Gaussian field $\mfX$ over a Hilbert space $\bH$ is a continuous linear
mapping $f\mapsto \mfX(f)$ from $\bH$ to the space of Gaussian random variables.
 The corresponding
chaos space $\HX$ is the Hilbert space of
square integrable random variables that are measurable with respect to the
sigma-algebra generated by $\mfX(f),\ f\in \bH$. The  chaos expansion
is an orthogonal decomposition of $\HX$: given an orthonormal
basis $\{\xi_m,\, m\geq 1\}$ in $\HX$, a square integrable
$\bH$-valued random variable $\eta$ has a chaos expansion
$\eta=\sum_{m\geq 1} \eta_m\xi_m$, with $\eta_m=\bE(\eta\xi_m)\in \bH$.

The definition of a generalized Gaussian field $\mfX$ already provides the
stochastic integral $\mfX(f)$ for non-random $f\in \bH$. As a result, given the
chaos expansion of a random element $\eta$ from $\HX$,
the definition of the stochastic integral $\mfX(\eta)$
requires an extension of the linearity property of $\mfX$ to linear
combinations with random coefficients. Two ``natural'' extensions of this
property lead to the \Ito-Skorokhod and the Stratonovich stochastic integrals;
see Definition \ref{def:main} below. Both integrals can be expressed using
the Malliavin derivative and divergence operator on the chaos space $\HX$.

Even for non-random $f$, when  there is no difference between the
\Ito-Skorokhod and the Stratonovich interpretations of $\mfX(f)$, there are
often several ways of computing $\mfX(f)$.
It is most convenient to work with a white noise over
 $\bH$, that is, a zero-mean generalize Gaussian field
such that $\bE\big(\mfX(f)\mfX(g)\big)=(f,g)_{\bH}$ for all $f,g\in \bH$.
It turns out that, for  every zero-mean Gaussian field $\mfX$ over $\bH$,
 there exists a different (usually larger) Hilbert space $\bH'$ such that $\mfX$ is a
white noise over $\bH'$.
Moreover, the space $\bH'$ is uniquely determined by $\mfX$.
On the other hand, every
 zero-mean Gaussian field $\mfX$ over $\bH$ can be written in the form
 $\mfX(f)=\mfB(\cK^*f)$, $f\in \bH$,
 where $\cK^*$ is a bounded linear operator on $\bH$ and $\mfB$ is a white noise
 over $\bH$, although this white noise
 representation of $\mfX$ is not necessarily unique.
 Thus, different white noise representations of $\mfX$ lead to
 different formulas for computing $\mfX(f),$ and the
 chaos expansion is an efficient way for deriving those formulas.
 In particular, for both deterministic and random $f$, chaos expansion
 provides an explicit formula for $\mfX(f)$ in terms of the Fourier
 coefficients of the integrand $f$.

To define stochastic integral with respect to a
 Gaussian process $X=X(t)$, $t\in [0,T]$,
we construct a Hilbert space $\bH_X$ and a white noise $\mfB$ over
$\bH_X$ such that
 $X(t)=\mfB(\chi_t)$, where $\chi_t$ is the characteristic
function of the interval $[0,t]$. The space $\bH_X$ is
 uniquely determined by $X$;
for example, the Wiener process on $(0,T)$ has $\bH_X=L_2((0,T))$.
 Then the equality
\bel{eq:def1}
\int_0^Tf(s)dX(s)=\mfB(f),\ f\in \HB,
\ee
  is a canonical
definition of the stochastic integral with respect to $X$.

  In some situations, given a Gaussian process $X=X(t),\ t\in [0,T]$,
  it is possible to find a
generalized Gaussian field $\mfX$ over a Hilbert space
$\bH$ so that $X(t)=\mfX(\chi_t)$. Even though  $\mfX$ is not
necessarily a white noise over $\bH$, the resulting definition
of the stochastic integral,
$$
\int_0^Tf(t)dX(t)=\mfX(f),
$$
 coincides with the (\ref{eq:def1}), while
the space $\bH$ can be more
convenient for computations than the space $\bH_X$. For example,
fractional Brownian motion  with the Hurst parameter bigger then $1/2$
has a rather complicated space  $\bH_X$, but
can be represented using a generalized Gaussian field over $\bH=L_2((0,T)).$

The paper is organized as follows.
Section 2 provides the definition and properties of generalized
Gaussian fields and establishes connections with the Gaussian processes.
Section 3 introduces the chaos expansion and the Wick product, both
necessary for the definition and analysis in Section 4 of the stochastic integrals
with random integrands.

The main contributions of the paper are:
\begin{enumerate}
\item Two white noise representations of a
zero-mean generalized Gaussian field (Theorem \ref{th:wn});
\item  A connection between  ge\-ne\-ra\-lized Gaussian
fields over $L_2((0,T))$ and processes that are representable in the
form $\int_0^t K(t,s)dW(s)$ (Theorem \ref{th:adapt});
\item Chaos expansions of the \Ito-Skorokhod and Stratonovich
integrals (Theorem \ref{th:st-in-ce});
\item Investigation of the equation $u(t)=1+\int_0^tu(s)dX(s)$ for
a class of Gaussian random processes $X$ (Theorem \ref{th:sode1}).
\end{enumerate}

 In particular, we establish the following result.

 \begin{theorem}
Let $\mfX$ be a zero-mean generalized Gaussian field over $L_2((0,T))$ and
$X(t)=\mfX(\chi_t)$. Then the solution of the \Ito~equation
$$
u(t)=1+\int_0^tu(s)dX(s)
$$ is unique in the class of square integrable $\cF^X$-measurable
processes and is given by
$$
u(t)=e^{X(t)-\frac{1}{2}\bE X^2(t)}.
$$
\end{theorem}

\section{Generalized Gaussian Fields}
\setcounter{theorem}{0}

Let $(\Omega,\,\cF,\, \bP)$ be a probability space and $\bA$, a
linear topological space over the real numbers $\bR$. Everywhere in this
paper, we assume that
the probability space is rich enough to support all the random
elements we might need.

\begin{definition}\label{Ch1:def1}
(a)  A \und{generalized random field} over $\bA$ is a
mapping \\
$\mfX: \Omega\times \bA \to \bR$ with the following properties:
\begin{enumerate}
\item For every $f\in \bA$, $\mfX(f)=\mfX(\cdot,f)$ is a random variable;
\item For every $\alpha, \beta\in \bR$ and $f,g\in \bA$,
$\mfX(\alpha f + \beta g)=\alpha \mfX(f)+\beta \mfX(g)$;
\item If \ $\lim\limits_{n\to \infty} f_n=f$ in the topology of $\bA$, then
$\lim\limits_{n\to \infty} \mfX(f_n)=\mfX(f)$ in probability.
\end{enumerate}

(b)  A generalized random field $\mfX$ is called
\begin{itemize}
\item \und{zero-mean}, if $\ \bE \mfX(f)=0$ for all $f\in \bA$;
\item  \und{Gaussian},
 if the random variable $\mfX(f)$ is Gaussian for
every $f\in \bA$.
\end{itemize}
\end{definition}

{\sc For Example}, if $W=W(t),\, 0\leq t\leq T,$ is a standard Brownian motion on
$(\Omega, \, \cF,\, \bP)$, then $\mfX(f)=\int_0^Tf(t)dW(t)$ is
a zero-mean generalized Gaussian field over $L_2((0,T))$; note that
\bel{ch1.1}
\bE|\mfX(f_n)-\mfX(f)|^2=\int_0^T|f_n(t)-f(t)|^2dt.
\ee
More generally, if ${\mathcal{M}}$ is  a bounded linear operator on $L_2((0,T))$,
then
\bel{ch1.2}
\mfX(f)=\int_0^T({\mathcal{M}}f)(t)dW(t)
\ee
is a zero-mean generalized Gaussian field over $L_2((0,T))$. In fact, by
 Theorem \ref{th:wn11}(b) below,
every zero-mean generalized Gaussian field over $L_2((0,T))$ can be represented in the
form (\ref{ch1.2}) with suitable ${\mathcal{M}}$ and $W$. We will also see
that the fractional Brownian motion on $[0,T]$
with Hurst parameter bigger than $1/2$
can be interpreted as a zero-mean generalized Gaussian field over $L_2((0,T)).$

 Let $\bH$ be a real Hilbert space with inner product $(\cdot, \cdot)_{\bH}$ and
 norm $\|\cdot\|_{\bH}=\sqrt{(\cdot,\cdot)_{\bH}}$. The following is a useful
 property of generalized Gaussian fields over $\bH$.

\begin{theorem}\label{th:gf}
For every zero-mean generalized Gaussian field $\mfX$  over a Hilbert space $\bH$,
 there exists a unique bounded linear self-adjoint operator $\cR$ on $\bH$ such that
\bel{ch1.3}
\bE\big( \mfX(f)\mfX(g)\big)=(\cR f,g)_{\bH},\ f,g \in \bH.
\ee
\end{theorem}

\pr Equality (\ref{ch1.3}) implies that the operator $\cR$,
if exists, must be unique.
To establish existence of $\cR$,
denote by $\wt{\bH}$ the Hilbert space $L_2(\Omega,\,\cF,\bP)$
of square integrable random variables. By Definition
\ref{Ch1:def1}, the mapping $ f\mapsto \mfX(f)$ defines a continuous
linear operator from $\bH$ to $\wt{\bH}$ (recall that, for
Gaussian random variables, convergence in probability implies
mean-square convergence). Therefore, there exists a positive number
$C$ such that, for every $f\in \bH$, \bel{Ch1-100}
\|\mfX(f)\|_{\wt{\bH}}^2 = \bE|\mfX(f)|^2 \leq C \|f\|_{\bH}^2.
\ee Fix $f\in \bH$ and consider the linear functional $F$  on
$\bH$ defined by $F(g) = \bE\big( \mfX(f)\mfX(g)\big)$. By
(\ref{Ch1-100}), this functional is bounded:
$$
|F(g)|=|\bE\big( \mfX(f)\mfX(g)\big)|\leq
\sqrt{\bE|\mfX(f)|^2} \sqrt{\bE|\mfX(g)|^2}\leq  C\|f\|_{\bH}\,\|g\|_{\bH},
$$
 and therefore, by the Riesz
Representation Theorem, there exists a unique $h_f\in \bH$ such
that $F(g)=\bE\big(\mfX(f)\mfX(g)\big)=(h_f,g)_{\bH}$.
Define the operator $\cR$ by $\cR f = h_f$. By
construction, this operator is linear; it is bounded by
(\ref{Ch1-100}). A bounded linear operator
satisfying (\ref{ch1.3}) automatically satisfies $ (\cR
f,g)_{\bH}=(f,\cR g)_{\bH}$ and is therefore self-adjoint.
\endproof

\begin{definition} (a) The operator $\cR$ from Theorem \ref{th:gf}
is called the \und{covariance operator} of $\mfX$.
(b)  A \und{white noise} over $\bH$ is a zero-mean generalized Gaussian field with
the covariance operator equal to the identity operator.
\end{definition}

Note that if $\cR$ is the covariance operator of $\mfX$ and $\cR f=0$,
  then $\mfX(f)=0$ ($\bP$-a.s.) Writing $\ker(\cR)$ to denote the zero-space
  of $\cR$, we have a direct sum decomposition $\bH=\ker(\cR)\bigoplus
  \ker(\cR)^{\perp}$, where $\ker(\cR)^{\perp}$ is the orthogonal
  complement of $\ker(\cR)$. As a result, if $f=f_1+f_2$, with
  $f_1\in \ker(\cR)$, $f_2\in \ker(\cR)^{\perp}$, then
  $\mfX(f)=\mfX(f_2)$. We say that the Gaussian field is \und{non-degenerate}
  if $\ker(\cR)=0$.

We show next that every zero-mean Gaussian random field over a Hilbert
space can be reduced to a white noise in two different ways.

\begin{theorem}\label{th:wn}
(a) For every zero-mean generalized Gaussian field  $\mfX$ over a Hilbert space $\bH$,
there exist a bounded linear  operator $\cK$ on $\bH$
 and a
white noise $\mfB$ over $\bH$ so that $\cK\cK^*$ is the covariance operator of $\mfX$
and, for every $f\in \bH$,
\bel{wn11}
\mfX(f)=\mfB(\cK^* f);
\ee
as usual, $\cK^*$ denotes the adjoint of $\cK$.

(b) For every zero-mean non-degenerate
generalized Gaussian field $\mfX$ over a Hilbert space $\bH$,
there exists a Hilbert space $\bH_{\cR}$ such that
$\bH$ is continuously embedded into $\bH_{\cR}$ and
$\mfX$ extends to a white noise over $\bH_{\cR}$.
\end{theorem}

\pr
(a) By construction, the covariance operator $\cR$ of a generalized  Gaussian field
 is non-negative definite, bounded, and
self-adjoint on $\bH$. Indeed, $\cR$ is bounded on $\bH$ by
Theorem \ref{th:gf}, and, for every $f$ and $g$ from $\bH$, we have
\begin{equation*}
\begin{split}
&(\cR f,f)_{\bH}=\bE\big(\mfX(f)\big)^2\geq 0;\\
&(\cR f,g)_{\bH}=\bE\big(\mfX(f)\mfX(g)\big)=
\bE\big(\mfX(g)\mfX(f)\big)=(\cR g,f)_{\bH}=(f,\cR g)_{\bH}.
\end{split}
\end{equation*}
  Therefore, by a standard result from functional analysis
  (see, for example, \cite[page 923]{DSch2}) there exists
a bounded linear operator $\cK$ on $\bH$ such that
$\cR= \cK\cK^*$; this operator $\cK$ is not necessarily unique.

 Next, let $\cK^*(\bH)=\{\cK^*f,\, f\in \bH\}$ be  the range of $\cK^*$, which
  is a closed linear subspace of
 $\bH$.
 Denote by $\cK^*(\bH)^{\perp}$  the orthogonal complement of $\cK^*(\bH)$
 in $\bH$.
Then, for every $f\in \bH$, there exists a unique pair $(f_1, f_2),$
 with  $f_1\in \bH$,  $f_2\in\cK^*(\bH)^{\perp}$,  such that
 $f=\cK^*f_1+f_2$. This orthogonal decomposition of $f$ implies
 \bel{aux-wn0002}
 \|f\|_{\bH}^2=\|\cK^*f_{1}\|^2_{\bH}+\|f_2\|^2_{\bH},
 \ee
 and
 \bel{aux-wn0003}
 \bE\big(\mfX(f_1)\big)^2=(\cK\cK^*f_1,f_1)_{\bH}=\|\cK^*f_1\|_{\bH}^2.
 \ee
Define
$$
\wt{\mfX}(f)=\mfX(f_1).
$$
Then $\wt{\mfX}$ is a generalized Gaussian field over $\bH$: if
$\lim\limits_{n\to\infty}\|f_n-f\|^2_{\bH}=0$, then, by
(\ref{aux-wn0002}) and (\ref{aux-wn0003}),
$$
\lim_{n\to \infty} \bE\big(\wt{\mfX}(f)-\wt{\mfX}(f_n)\big)^2=
\lim_{n\to \infty} \bE\big(\mfX(f_1-f_{1,n})\big)^2=
\lim_{n\to \infty}\|\cK^*(f_1-f_{1,n})\|^2_{\bH}=0.
$$
 Let $\ol{\mfB}$ be a white noise over $\bH$,
  independent of $\mfX$. The same arguments show that $\wt{\mfB}$, defined by
  $$
  \wt{\mfB}(f)=\ol{\mfB}(f_2),
  $$
  is a generalized Gaussian field over $\bH$.
  Define
 $$
\mfB(f)=\wt{\mfX}(f)+\wt{\mfB}(f).
 $$
 Then $\mfB$ is a generalized Gaussian field over $\bH$,
 being a sum of two independent
 generalized Gaussian
 fields over $\bH$, and, by definition, $\mfB(\cK^*f)=\wt{\mfX}(\cK^*f)=\mfX(f)$.
  Moreover, if
  $f=\cK^* f_1+f_2$, $g=\cK^*g_1+g_2$, then
  \begin{equation*}
  \begin{split}
  \bE\big(\mfB(f)\mfB(g)\big)&=\bE\big(\mfX(f_1)\mfX(g_1)\big)+
  \bE\big(\ol{\mfB}(f_2)\ol{\mfB}(g_2)\big)\\ &=
  (\cK\cK^* f_1,g_1)_{\bH}+(f_2,g_2)_{\bH}
  = (\cK^*f_1,\cK^*g_1)_{\bH}+(f_2,g_2)_{\bH}=
  (f,g)_{\bH},
  \end{split}
  \end{equation*}
  where the first equality follows from the independence of $\mfX$ and
  $\wt{\mfB}$, and the last, from $(\cK^*f_1,g_2)_{\bH}=(\cK^*g_1,f_2)_{\bH}=0$.
 Thus, $\mfB$ is a white noise over $\bH$, and the proof of
 (\ref{wn11}) is complete.

 (b) Define $\bH_{\cR}$ as the closure of $\bH$
 with respect to the inner product $(f,g)_{\bH_{\cR}}=(\cR f,g)_{\bH}$.
 For  $f\in \bH$,
 $\|f\|_{\bH_{\cR}}^2=(\cR f,f)_{\bH}\leq C\|f\|_{\bH}^2$, which
  implies a
 dense continuous embedding of $\bH$ into $\bH_{\cR}$.  By definition,  for
 $f,g\in \bH$, $\bE\big(\mfX(f)\mfX(g)\big)=(\cR f,g)_{\bH}=(f,g)_{\bH_{\cR}}$.
 As a result, if $f\in \bH_{\cR}$ and $\lim_{n\to \infty}\|f_n-f\|_{\bH_{\cR}}^2=0$,
 with $f_n\in \bH$, then
 $$
 \lim_{m,n\to \infty}\bE\big(\mfX(f_m)-\mfX(f_n)\big)^2=
\lim_{m,n\to \infty} \|f_m-f_n\|_{\bH_{\cR}}^2=0,
$$
so that $\lim\limits_{n\to \infty}\mfX(f_n)$ exists in the mean-square
and is therefore a Gaussian random variable.
 We then define $\mfX(f)=\lim\limits_{n\to \infty}\mfX(f_n)$.
 The value of $\mfX(f)$ does not depend on the
sequence $\{f_n,\, n\geq 1\}$ approximating $f$, because
$$
\lim_{n\to \infty}\bE\big(\mfX(f)-\mfX(f_n)\big)^2 = \|f-f_n\|_{\bH_{\cR}}^2=0.
$$
Also,
$$
\bE\big(\mfX(f)\mfX(g)\big)=\lim_{n\to \infty}
\bE\big(\mfX(f_n)\mfX(g_n)\big)=\lim_{n\to \infty}( f_n,g_n)_{\bH_{\cR}}=
(f,g)_{\bH_{\cR}},
$$
meaning that this extension of $\mfX$ is a white noise over $\bH_{\cR}$.
\endproof

\begin{remark}\label{rem-wn}
 (a) If $\mfX$ is non-degenerate and $\cR: \bH\to\bH$
 is onto, then $\cR$ has a bounded inverse and
  $\bH_{\cR}=\bH$.
  (b) If $\ker{\cR}$ is non-trivial, then
we can define $\bH_{\cR}$ as the closure of the factor space $\bH/\ker(\cR)$
 with respect to the inner product $(\ol{f},\ol{g})_{\bH_{\cR}}=(\cR f,g)_{\bH}$,
 where $\ol{f}$ is the equivalence class of $f$ in $\bH/\ker(\cR)$.
 Direct computations show that
 the  generalized random field $\mfB$ over $\bH/\ker(\cR)$, defined by
 $$
 \mfB(\ol{f})=\mfX(f), \ f\in \bH,
 $$
 extends to a white noise over $\bH_{\cR}$.
 \end{remark}

We will now discuss several connections between generalized
Gaussian fields  and Gaussian processes. In what follows, $I$ denotes
 either an interval $[0,T]$, or the half-line
$[0,+\infty)$, or all of $\bR$.

 Denote by $\chi_t=\chi_t(s)$  the
characteristic function of the interval $[0,t]$:
\bel{char-f}
\chi_t(s)=
\begin{cases}
1, & 0\leq s \leq t;\\
0, & {\rm otherwise}.
\end{cases}
\ee
With this definition, $\chi_{t_2}(s)-\chi_{t_1}(s)$ is the
characteristic function of the interval $(t_1,t_2]$, $t_2>t_1$.

\begin{theorem} \label{th:wn11}
 (a) If $\mfB$ is a white noise over $L_2(I)$,
 then $B(t)=\mfB(\chi_t)$ is a standard Brownian motion on
$(\Omega, \, \cF,\, \bP)$ and, for every $f\in L_2(I)$, we have
\bel{wn-rep1}
\mfB(f)=\int_If(s)dB(s).
\ee
(b) For every zero-mean non-degenerate generalized Gaussian field
$\mfX$ over $L_2(I)$, there exist a bounded linear operator $\cK^*$ on $L_2(I)$ and a
standard Brownian motion $W=W(t)$ such that, for every $f\in L_2(I)$,
\bel{wn-rep2}
\mfX(f)=\int_I(\cK^*f)(s)dW(s).
\ee
\end{theorem}

\pr (a) Direct computations show that $B=B(t)$ has all the properties of the
standard Brownian motion. In particular,
$$
\bE\big( B(t_1)B(t_2)\big)=\int_0^T\chi_{t_1}(s)\chi_{t_2}(s)ds=\min(t_1,t_2).
$$
Next, if  $s_0 < s_1< \ldots< s_N$ is a finite collection of points in $I$ and
$f(s)=\sum\limits_{k=1}^N a_k(\chi_{s_k}(s)-\chi_{s_{k-1}}(s))$ is a
(non-random) step function, then the linearity property of the generalized
random field $\mfB$ implies
$$
\mfB(f)=\sum_{k=1}^Na_k\big(\mfB(\chi_{s_k})-\mfB(\chi_{s_{k-1}})\big)=\int_If(s)dB(s).
$$
For general $f$, the result then follows after passing to the limit,
 using the continuity property of the
 generalized random field $\mfB$ and the $L_2$-isometry of the
 stochastic integral.

 (b) This follows from part (a) and from Theorem \ref{th:wn}.
 \endproof

 Given a zero-mean generalized Gaussian field $\mfX$ over $L_2(I)$, we define its
 \und{associated process} $X(t),\ t\in I,$ by
 \bel{a-proc}
 X(t)=\mfX(\chi_t).
 \ee
 Clearly, $X(t)$ is a Gaussian process. Let $\cK^*$ be the operator
 from Theorem \ref{th:wn11} and define the kernel function
 $K_{\mfX}=K_{\mfX}(t,s)$ by
 \bel{functionK}
K_{\mfX}(t,s)=(\cK^*\chi_t)(s).
\ee
 It then follows from (\ref{wn-rep2}) that
 \bel{a-proc1}
 X(t)=\int_I K_{\mfX}(t,s)dW(s)
 \ee
 for some standard Brownian motion $W$. Let us emphasize that, while
 every kernel $K(t,s)$ with minimal integrability properties can
 define a Gaussian process according to (\ref{a-proc1}), only a process
 associated with a generalized field over $L_2(I)$ has
 a kernel defined according to (\ref{functionK}), where $\cK^*$ is a
 bounded operator on $L_2(I)$. Recall that the definition
 of a generalized field (Definition \ref{Ch1:def1})
 includes a certain continuity property, and this property translates into addition
 structure of the kernel function in the representation of the associated
process.

Now assume that we are given a Gaussian process $X(t)$ defined by
(\ref{a-proc1}) with some kernel $K_{\mfX}(t,s)$. {\em We are not assuming
that $K_{\mfX}$ has the form (\ref{functionK})}.  In what follows, we
discuss sufficient conditions on $K_{\mfX}(t,s)$ ensuring that
$X(t)$
is the associated process of a generalized Gaussian field $\mfX$ over
$L_2(I)$, that is,
representation (\ref{functionK}) does indeed hold with some bounded linear operator
$\cK^*$ on $L_2(I)$.
For that, we need to  recover the operator $\cK^*$
from the kernel $K_{\mfX}(t,s)$.  By linearity, if  (\ref{functionK}) holds and
 if $s_0<s_1<\ldots<s_N$ are points in $I$ and
\bel{stepF}
f(s)=\sum\limits_{k=0}^{N-1} a_k(\chi_{s_{k+1}}(s)-\chi_{s_{k}}(s))
\ee is a step function, then
\bel{op1}
\cK^* f(s)= \sum_{k=0}^{N-1} a_k \big(K_{\mfX}(s_{k+1},s)-K_{\mfX}(s_k,s)\big).
\ee
To extend (\ref{op1}) to continuous functions $f$, the kernel
$K_{\mfX}(t,s)$ must have bounded variation as a function of $t$;
if this is indeed the case, then (\ref{op1}) implies that, for every
smooth compactly supported function $f$ on $I$,
\bel{op2}
\cK^* f(s)=\int_I f(t)K_{\mfX}(dt,s).
\ee
The assumption about the bounded variation of the kernel is used
extensively in \cite{AMN}, and the connection with generalized fields
 shows that this assumption is very natural. It now follows that if
 the partial derivative $\partial K_{\mfX}(t,s)/\partial t$ exists and
 is square integrable over $I\times I$, then $\cK^*$,
 as defined by (\ref{op2}), extends to a bounded linear operator on $L_2(I)$.

 Let us now assume that $I=[0,T]$ and the process $X(t)$ define by (\ref{a-proc1})
 is \und{non-anticipating},  i.e.
 adapted to the filtration $\{\cF_t^W,\ 0\leq t\leq T\}$
  generated by the Brownian motion $W(s).$  Then $K_{\mfX}(t,s)=0$ for $s>t$ and
 (\ref{a-proc1}) becomes
 \bel{a-proc2}
 X(t)=\int_0^t K_{\mfX}(t,s)dW(s).
 \ee
 Note that in this case we have
 \bel{cov-funct}
 \bE\big(X(t)X(s)\big)=\int_0^{\min(t,s)}K_{\mfX}(t,\tau)K_{\mfX}(s,\tau) d\tau.
 \ee
 This is the type of processes studied in \cite{AMN}, and for such processes,
 formula (\ref{op1}) and the conditions for the continuity of the
 corresponding operator $\cK^*$ must be modified as follows.

 \begin{theorem}\label{th:adapt}
 Assume that $I=[0,T]$ and
 the process $X(t)$ defined by (\ref{a-proc1}) is non-anticipating.

 (a) If $f$ is a step function (\ref{stepF}), then
 \bel{KK}
\begin{split}
\cK^*f(s)&=\sum_{i=0}^{N-1}\big(\chi_{s_{i+1}}(s)-\chi_{s_i}(s)\big)
\Big(a_iK_{\mfX}(s_{i+1},s)\\
&+\sum_{k=i+1}^{N-1}a_{k}\big(K_{\mfX}(s_{k+1},s)-
K_{\mfX}(s_{k},s)\big)
\Big).
\end{split}
\ee
(b) If the function $K_{\mfX}(\cdot, s)$ has bounded
variation for every $s$ and  $\lim\limits_{\delta\to 0,\,
\delta>0}K_{\mfX}(s+\delta,s)=K_{\mfX}(s^+,s)$ exists for all $s\in (0,T),$
 then
  \bel{ker-op000}
   \cK^*f(s)=K_{\mfX}(s^+,s)f(s)+\int_s^T f(t)K_{\mfX}(dt,s)
 \ee
 for every continuous on $[0,T]$ function $f$.

 (c) If the function
$K_{\mfX}(t,s)$ has the following properties
\begin{enumerate}
\item $K$ is continuous and non-negative for  $0\leq s\leq t \leq T$,
 and $\sup\limits_{0<t<T} K(t,t) \leq K_0$;
\item $K^{(1)}(t,s)=\partial K(t,s)/\partial  t$ is non-negative for $0< s < t < T$
and there exists a number $K_1=K_1(T)$ such that
\bel{ker-cond1}
\sup_{0<t<T} \int_0^t K(T,s)  K^{(1)}(t,s) ds \leq K_1(T),
\ee
\end{enumerate}
then the corresponding operator $\cK^*$ defined by equation
\eqref{ker-op000} is bounded on $L_2((0,T))$ and the operator
norm $\|\cK^*\|$ of $\cK^*$ satisfies
\bel{op-norm}
\|\cK^*\|^2 \leq
\begin{cases}
2(K_0^2+K_1),& {\rm if \ } K_0>0;\\
K_1,& {\rm if \ } K_0=0.
\end{cases}
\ee
\end{theorem}

\pr (a) By assumption, $K_{\mfX}(t,s)=0$ for $s>t$.
 Fix an $s$ such that $s\in (s_j,s_{j+1}]$ for some $j=0,\ldots, N-1$.
  By (\ref{op1}) we have for this value of $s$
\begin{equation*}
\begin{split}
\cK^*f(s)&=\sum_{k=0}^{N-1} a_k \big(K_{\mfX}(s_{k+1},s)-
K_{\mfX}(s_{k},s) \big)\\
&=a_{j}K_{\mfX}(s_{j+1},s)
+\sum_{k=j+1}^{N-1}a_{k}\big(K_{\mfX}(s_{k+1},s)-
K_{\mfX}(s_{k},s)\big).
\end{split}
\end{equation*}
Since $\chi_{s_{k+1}}(s)-\chi_{s_{k}}(s)$ is the characteristic function
of the interval $(s_k,s_{k+1}]$, (\ref{KK}) follows.

(b) Under the additional assumptions on the kernel $K_{\mfX}$,
 (\ref{ker-op000}) follows from (\ref{KK}) after passing to the limit
 $\max\limits_{j=0,\ldots, N-1}|s_{j+1}-s_j|\to 0$.

 (c) Let $g$ be a smooth compactly supported function on $(0,T)$.
 It follows from (\ref{ker-op000}) that
 $$
 \cK^*g(s)=
 K(s,s)g(s)ds+\!\int_s^T\!\frac{\partial K(\tau,s)}{\partial  \tau}\, g(\tau)\, d\tau=
K(s,s)g(s)ds+\!\int_s^T \!K^{(1)}(\tau,s)\, g(\tau)\, d\tau.
 $$
  To estimate the $L_2$-norm of the integral, we use the Cauchy-Schwartz
  inequality and the properties of $K^{(1)}$:
  \begin{equation*}
\begin{split}
&\int_0^T \left| \int_s^T  K^{(1)}(\tau,s)
                g(\tau) d\tau \right|^2 ~ds
                = \int_0^T \left| \int_s^T \left[K^{(1)}(\tau,s)\right]^{1/2}
                \left[K^{(1)}(\tau,s)\right]^{1/2} g(\tau) d\tau \right|^2 ~ds \\
            &\leq \int_0^T \int_s^T K^{(1)}(\tau,s) d\tau
                \int_s^T K^{(1)}(\tau,s) g^2(\tau) d\tau ~ds\\
           & \leq \int_0^T \left( K(T,s) - K(s,s)\right) \int_s^T
                K^{(1)}(\tau,s) g^2(\tau) d\tau ~ds \notag\\
            &\leq \int_0^T \left( \int_0^\tau K(T,s)
                K^{(1)}(\tau,s) ~ds \right) g^2(\tau) ~d\tau \notag
    \leq K_1(T) \|g\|_{L_2((0,T))}^2.
    \end{split}
    \end{equation*}
    \endproof

We remark that in \cite{AMN} relation (\ref{KK}) is used to {\em define} the
operator $\cK^*$ corresponding to a non-anticipating process $X(t)$.
 Using the connection with the generalized fields, Theorem \ref{th:adapt}
shows that this definition is reasonable.

    The main example covered by part (c) of Theorem \ref{th:adapt} is
the fractional Brownian motion $W^H$ on $[0,T]$ with the Hurst parameter $H>1/2$.
Indeed, it is known (see \cite[Section 5.1.3]{Nualart})  that in this case $W^H$
has representation
(\ref{a-proc2}) with
$$
K_{\mfX}(t,s)=C_H\left(H-\frac{1}{2}\right)
 s^{\frac{1}{2}-H}\int_s^t (\tau-s)^{H-\frac{3}{2}}\tau^{H-\frac{1}{2}}\, d\tau,
$$
where
$$
C_H=\left(\frac{2H\Gamma\left(\frac{3}{2}-H\right)}{\Gamma\left(H+\frac{1}{2}\right)
\Gamma(2-2H)}\right)^{\frac{1}{2}}
$$
and $\Gamma$ is the Gamma-function. Clearly, $K_{\mfX}(s,s)=0$ and so $K_0=0$.
 Then somewhat lengthy computations show that
\bel{fbm-bound}
K_1(T)=\frac{H(2H-1)\,\Gamma\left(H-\frac{1}{2}\right)}
{\Gamma\left(H+\frac{1}{2}\right)}\ T^{2H-1}.
\ee
The bound  $K_1(T)$ is asymptotically optimal:
 since $\lim\limits_{x\to 0^+}
x\Gamma(x)=\lim\limits_{x\to 0^+}\Gamma(1+x)=1$,
the right-hand side of (\ref{fbm-bound}) converges to $1$ as
 $H\searrow \frac{1}{2}$, and
if $H=1/2$, then $W^H$ is the standard Brownian
motion and $\cK^*$ is the identity operator, which corresponds to $\|\cK^*\|=1$.

 The following theorem establishes a connection between a zero-mean Gaussian
 process and white noise.

 \begin{theorem} \label{th:gp-wn} For every zero-mean  Gaussian process
 $X=X(t), \, t\in I\subseteq \bR $ with covariance
 function $R(t,s)=\bE\big(X(t)X(s)\big)$, there exist
 \begin{enumerate}
 \item
 a Hilbert space $\bH_{R}$ containing the indicator functions $\chi_t;$
 \item  a white noise $\mfB$ over $\bH_{R}$
 \end{enumerate}
 such that $X(t)=\mfB(\chi_t)$.
 \end{theorem}

 \pr Let the Hilbert space
 $\bH_{R}$  be the closure of the set of the step functions with respect to
 the inner product
 $$
 (\chi_{t_1},\chi_{t_2})_{\bH_{R}}=R(t_1,t_2).
 $$
 Define a generalized
 Gaussian field $\mfB$ over $\bH_{R}$ by setting
 \bel{gen-f-p}
 \mfB(\chi_t)=X(t),
 \ee
   and then extending by linearity and
  continuity to all of $\bH_{R}$. With this definition, $\mfB$ is
 a white noise over $\bH_R$.
\endproof

By analogy with (\ref{wn-rep1}), if $\mfX$ is a generalized Gaussian field
over a Hilbert space $\bH$ of functions or generalized functions on $I$, and
 $X(t)$ is the associated process
of $\mfX$, then $\int_I f(s)dX(s)$ can be an alternative notation for $\mfX(f)$.

 The space $\bH_{R}$ from Theorem \ref{th:gp-wn} appears in \cite{AMN} and
  is different from
 reproducing kernel Hilbert space used in \cite[Section 6]{PT}.
 If $X(t)$ is the associated process of a zero-mean non-degenerate
 generalized Gaussian field
 $\mfX$ over $\bH=L_2(I)$, and $\cR$ is the covariance operator of $\mfX$, then
 $R(t,s)=(\cR \chi_t,\chi_s)_{L_2(I)}$ and the space $\bH_{R}$
 coincides with $\bH_{\cR}$ from Theorem \ref{th:wn}.

 Given a covariance function $R$, an explicit characterization of the
  space $\bH_R$ is impossible without additional assumptions about
  $R$. For example, in \cite{AMN},  representation
  $$
  R(t,s)=\int_0^{\min(t,s)}K(t,\tau)K(s,\tau)d\tau,
  $$
   is used, along with various assumptions about the kernel $K$.
 If $I=[0,T]$ and $R(t,s)=\min(t,s)$, then
 $(\chi_{t_1},\chi_{t_2})_{\bH_{R}}=(\chi_{t_1},\chi_{t_2})_{L_2((0,T))}$.
 That is, for the Wiener process,  $\bH_R=L_2((0,T))$.

  Let us summarize the main results of this section:
  \begin{itemize}
  \item Every zero-mean generalized Gaussian random field over $\bH$ with covariance operator
  $\cR$ has two  white noise representations:
  over the Hilbert space $\bH_{\cR}$ and  over the original space $\bH$;
  \item Every zero-mean Gaussian random process with covariance function $R$ is
  the associated process of a white noise over the Hilbert space $\bH_{R}$.
  \end{itemize}

\section{Chaos Decomposition and the Wick Product}
\setcounter{theorem}{0}

Let $\mfX$ be a zero-mean generalized Gaussian field over a real
Hilbert space $\bH$, on a probability space $(\Omega,\, \cF,\, \bP)$. From now
on, we assume that the space $\bH$ is separable. Denote by
$\cF^{\mfX}$ the sigma-algebra generated by the random variables
$\mfX(f),\ f\in \bH$.

\begin{definition} (a) The \und{chaos space} generated by $\mfX$ is the
collection of all random variables on $(\Omega,\, \cF,\, \bP)$ that
are square integrable and $\cF^{\mfX}$-measurable. This chaos space
will be denoted by $\HX$.

(b) The \und{first chaos space} generated by $\mfX$ is the sub-space of $\HX$,
consisting of the random variables $\mfX(f)$, $f\in \bH$.
 The first chaos space will be denoted by $\HXt$.
\end{definition}

It follows that $\HX$ is a Hilbert space with inner product $(\xi,\eta)_{\HX}=
\bE(\xi\eta)$, and $\HXt$ is a Hilbert sub-space of $\HX$.
Moreover, the space $\HXt$ is separable: if
$\{\bar{f}_1, \bar{f}_2, \ldots\}$ is
 a dense countable set in $\bH$, then the collection of all finite linear combinations
of $\mfX(\bar{f}_i)$ with rational coefficients is a dense countable set in
$\HXt$.

Our next objective is to show how an orthonormal basis in $\HXt$ leads to an
orthonormal basis in $\HX$. We will need some additional constructions.

For an integer $n\geq 0$,
  the $n$-th \und{Hermite polynomial} $H_n=H_n(t)$ is defined by
\begin{equation}
\label{eq:HerPol}
 H_{n}(t)=(-1)^ne^{t^{2}/2}
\frac{d^{n}}{dt^{n}}e^{-t^{2}/2}.
\end{equation}
 In particular,   $H_0(t)=1$, $H_1(t)=t$, $H_2(t)=t^2-1$, $H_3(t)=t^3-3t$, etc.
 Note that $H_n(t)=t^n+\ldots$, that is, $H_n$ is a polynomial
 of degree $n$ and  the leading coefficient is always equal to one.
 It is well known that if $\xi$ is a standard Gaussian random variable, then
 \bel{orthon}
 \bE\big( H_n(\xi)H_m(\xi) \big) =
\begin{cases} n!,& n=m;\\
0,& n\not=m.
\end{cases}
\ee
In fact, the collection $\{ H_n(\xi),\ n\geq 0\}$ is an orthonormal
basis in the space of square integrable, $\cF^{\xi}$-measurable random
variables.

Next, denote by $\cI$ the collection of
multi-indices, that is, sequences
 $\alpha=\{\alpha_k,\ k\geq 1\}= \{\alpha_1, \alpha_2,\ldots\}$
 with the following properties:
 \begin{itemize}
 \item each $\alpha_k$ is a non-negative integer: $\alpha_k\in \{0,1,2,\ldots\}$.
 \item only finitely many of $\alpha_k$ are non-zero:
 $|\alpha|:=\sum\limits_{k=1}^{\infty} \alpha_k < \infty. $
 \end{itemize}
 The set $\cI$ is countable, being a countable union of countable sets.
  By $\epsilon_n$ we denote the multi-index $\alpha=\{\alpha_k,\ k\geq 1\}$
 with $\alpha_k=1$ if $n=k$ and $\alpha_k=0$ otherwise.
 For $\alpha \in \cJ$, we will use the notation
 $$
 \alpha!:=\alpha_1!\,\alpha_2!\cdots
 $$

 Let $\{\xi_1,\xi_2,\ldots\}$ be an ordered countable collection of random variables.
 For $\alpha\in \cI$ define random variables $\xi_{\alpha}$ as follows:
 \bel{xi-al}
 \xi_{\alpha}=\prod_{k\geq 1}\frac{H_{\alpha_k}(\xi_k)}{\sqrt{\alpha_k!}},
 \ee
 where $H_{\alpha_k}$ is $\alpha_k$-th Hermite polynomial (\ref{eq:HerPol}).
 For example,   $\alpha\!=\!(0,2,0,1,3,0,0,\ldots)$ has   three
non-zero entries $\alpha_2=2$, $\alpha_4=1$, and $\alpha_5=3$, so that
$$
\xi_{\alpha}=\frac{H_2(\xi_2)}{\sqrt{2!}}\cdot  H_1(\xi_4) \cdot
\frac{H_3(\xi_5)}{\sqrt{3!}}=\frac{\xi_2^2-1}{\sqrt{2}}\,\xi_4\,
\frac{\xi_5^3-3\xi_5}{\sqrt{6}}.
$$
 The product on the right  hand side of (\ref{xi-al})
 is finite for every $\alpha\in \cI$.
 Note also that $\xi_k=H_1(\xi_k)=\xi_{\epsilon_k}$ and, more generally,
  $H_n(\xi_k)= \sqrt{n!}\,\xi_{n\epsilon_k}$.

 The following theorem has been known for some time in various
 forms.
 In the particular case when $\mfX(f)=\int_0^T f(t)dW(t)$, this theorem
is the main result of the paper \cite{CM} by Cameron and Martin;
see also \cite[Theorem 1.9]{HKPS} and \cite[Theorem 2.2.3]{HOUZ}.
The  formulation and proof below are similar to \cite[Theorem 2.1]{Maj}.

 \begin{theorem}\label{th:basis}
 Let $\{\xi_1,\xi_2,\ldots\}$ be an orthonormal basis in $\HXt$.
 Then the collection $\Xi=\{\xi_{\alpha}, \ \alpha \in \cI\}$
 is an orthonormal basis in $\HX$: for every $\eta\in \HX$ we have
 $$
 \eta=\sum_{\alpha \in\cI} \Big(\bE(\eta\xi_{\alpha})\Big)\, \xi_{\alpha},\ \
 \bE\eta^2=\sum_{\alpha\in \cI} \Big(\bE(\eta\xi_{\alpha})\Big)^2.
 $$
 \end{theorem}

 \pr
Recall that $\mfX$ is a Gaussian random field. As a result,
 an orthonormal basis in $\HXt$ is a collection of standard
 Gaussian random variables $\xi_k$, $k\geq 1$,  that are uncorrelated,
 hence independent. Then property (\ref{orthon}) of Hermite polynomials
 implies that $\Xi$ is an orthonormal system.

 Next, denote by $H_{\xi_k}$ the Hilbert space of
 square integrable random variables that are measurable
 with respect to the sigma-algebra generated by $\xi_k$.
 Consider the product space $H_{\infty}=\prod\limits_{k=1}^{\infty} H_{\xi_k}$.
 By the definition of the product topology, it follows that the
 collection $\{\xi_{\alpha}, \ \alpha\in \cI\}$ is an orthonormal
 basis in this product space. We also
 note that, by construction, the sigma-algebra $\cF^{\mfX}$ is generated by
 the random variables $\xi_1, \xi_2,\ldots,$ and therefore, for every
 $\cF^{\mfX}$-measurable random variable $\eta$, there exists a measurable,
 real-valued function $F$ on the measurable space
 $(\bR^{\infty}, \mathcal{B}(\bR^{\infty}))$ with the property $\eta=F(\xi_1,\xi_2,
 \ldots)$.
 This establishes a one-to-one correspondence between the product space
 $H_{\infty}$ and the chaos space $\HX$, and completes the proof.
\endproof

By definition, the space $\HXt$ is generated by $\xi_{\alpha}$ with
$|\alpha|=1$. More generally, we define $\HX^{(N)}$,
 the \und{$N$-th chaos space} of $\mfX$,
  as the closure in $\HX$
  of the linear span of $\xi_{\alpha}$ with $|\alpha|=N$:
  $\eta\in\HX^{(N)}$ if and only if $\eta=\sum\limits_{\alpha \in \cI, |\alpha|=N}
  c_{\alpha} \xi_{\alpha}$ for some real numbers $c_{\alpha}$ satisfying
  $\sum_{\alpha} |c_{\alpha}|^2< \infty$.
 By Theorem \ref{th:basis} we have the  \und{chaos
 decomposition} of $\HX$:
 \bel{OrthDec}
 \HX=\bigoplus_{N=0}^{\infty}\HX^{(N)}=\HX^{(0)}\oplus\HX^{(1)}\oplus
 \HX^{(2)}\oplus\cdots.
 \ee

 \begin{proposition}
  For each $N$, the
   space $\HX^{(N)}$ does not depend on the choice of the basis in
   $\HXt$.
   \end{proposition}

 \pr
Since the polynomials $H_0, H_1,\ldots, H_N$ are orthogonal
with respect to the Gaussian measure on $\bR$, these polynomials
are linearly independent. Therefore, for each $N\geq0$, the space
$
\HX^{\leq N}=\HX^{(0)}\oplus\HX^{(1)}\oplus\cdots \oplus \HX^{(N)}
$
coincides with the closure in $\HX$ of the linear span of the random
variables
$P_N(\mfX(f_1),\ldots, \mfX(f_k))$, $k\geq 1$, $f_i\in \bH$,
where $P_N$ is a polynomial of degree at most $N$; cf. \cite[p.~9]{Maj}.
Thus, the space $\HX^{\leq N}$ does not depend on the basis in
$\HXt$. Since $\HX^{\leq N}=\HX^{\leq (N-1)}\oplus \HX^{(N)}$,
the space $\HX^{(N)}$ does not depend on the basis as well.
\endproof

In the case of white noise $\mfB$ over $\bH$, an orthonormal basis in $\HB^{(1)}$
is closely related to an orthonormal basis in $\bH$.

\begin{proposition} \label{prop:wne}
Let $\mfB$ be a white noise over a separable Hilbert space $\bH$ and let
$\{m_1,\, m_2,\, \ldots\}$ be an orthonormal basis in $\bH$. Then
$\{\xi_k=\mfB(m_k),\ k\ge 1\}$ is an orthonormal basis in $\HBt$ and, for
every $f\in \bH$,
\bel{wne}
\mfB(f)=\sum_{k=1}^{\infty} (f,m_k)_{\bH}\, \mfB(m_k).
\ee
\end{proposition}

\pr Note that $\bE \big(\xi_k \xi_n\big)=\bE \big(\mfB(m_k)\mfB(m_n)\big)
=(m_k,m_n)_{\bH}$, so
the system $\{\xi_k,\ k\ge 1\}$ is orthonormal in $\HBt$
 if and only if $\{m_k,\, k\geq 1\}$ is orthonormal in $\bH$.
If $\xi\in \HBt$, then $\xi=\mfB(f)$ for some $f\in \bH$.
By assumption, $f=\sum_{k=1}^{\infty}(f,m_k)_{\bH}\, m_k$, which implies
(\ref{wne}) and completes the proof. \endproof

If $\bH=L_2((0,T))$ and $f=\chi_t$, then (\ref{wne}) becomes
a familiar representation of the standard Brownian motion on $[0,T]$:
\bel{wn:wp}
W(t)=\sum_{k=1}^{\infty} \left(\int_0^t m_k(s)ds\right) \, \left(\int_0^T m_k(s)dW(s)
\right).
\ee
Now, let $\mfX$ be a zero-mean generalized Gaussian field over
a separable Hilbert space $\bH$.
By (\ref{wne}) and Theorem \ref{th:wn}(a),
we can take a white noise representation of $\mfX$,
$\mfX(f)=\mfB(\cK^* f)$, and get an expansion of $\mfX(f)$ using an
orthonormal basis in $\bH$:
\bel{wne1}
\mfX(f)=\sum_{k=1}^{\infty} (\cK^* f , m_k)_{\bH}\, \mfB(m_k).
\ee
When $\bH=L_2((0,T))$ and $f=\chi_t$, the associated process
has representation $X(t)=\int_0^TK_{\mfX}(t,s)dW(s)$,
where $K_{\mfX}(t,s)=(\cK^*\chi_t)(s)$,
 and we get a generalization of (\ref{wn:wp}):
\bel{asopr11}
X(t)=\sum_{k=1}^{\infty} \left(\int_0^t (\cK m_k)(s)ds\right) \,
\left(\int_0^T m_k(s)dW(s)\right);
\ee
note that $ \int_0^t(\cK m_k)(s)ds= \int_0^TK_{\mfX}(t,s)m_k(s)ds $.

Alternatively, by Theorem \ref{th:wn}(b) and Remark \ref{rem-wn}(b),
$\mfX$ is a white noise over the space $\bH_{\cR}$ corresponding to the
covariance operator $\cR$ of $\mfX$.
 If $\{\ol{m}_k, \, k\geq 1\}$
is an orthonormal basis in $\bH_{\cR}$, then
we have the following analog of (\ref{wne}):
\bel{nwne}
\mfX(f)=\sum_{k=1}^{\infty} (\cR f,\ol{m}_k)_{\bH}\, \mfX(\ol{m}_k).
\ee

If $\mfX$ is non-degenerate, which means $\ker(\cR)=0$,
then (\ref{wne1}) and (\ref{nwne}) are equivalent. Indeed,
 by Theorem \ref{th:wn}(b), $\bH$ is dense in $\bH_{\cR}$ and
 we can extend  $\cK^*$ to a bounded linear operator from $\bH_{\cR}$ to
 $\bH$, because $\|\cK^*f\|_{\bH}^2=\|f\|_{\bH_{\cR}}^2$. Clearly,
 (\ref{wne1}) and (\ref{nwne}) coincide for $f\in \bH$, since
 $\{\cK^* \ol{m}_k,k\geq 1\}$ is an orthonormal basis in $\bH$.
Extending  (\ref{wne1}) to $f\in \bH_{\cR}$ makes
(\ref{wne1}) equivalent to (\ref{nwne}).

We conclude the section with a brief discussion of the Wick product, as we will
need this product to define $\mfX(f)$ for random $f$.

To motivate the definition of the Wick product, we make the following observation.
The ordinary powers $x^n$ have the property $x^nx^m=x^{m+n}$. By Theorem
\ref{th:basis}, the
natural building blocks of the chaos space $\HX$ are not the ordinary powers
but Hermite polynomials of the basis elements in $\HXt$.
It is therefore convenient to have an operation,
which we denote by $\diamond$ and call the \und{Wick product}, so that, for
every $\xi\in \HXt$,
\bel{WPmot}
H_n(\xi)\diamond H_m(\xi)=H_{m+n}(\xi).
\ee
In fact, together with Theorem \ref{th:basis}, relation (\ref{WPmot})
completely defines the Wick product in $\HX$, because
if $\{\xi_{\alpha}, \, \alpha\in \cI\}$ is an orthonormal basis
in $\HX$, as defined by (\ref{xi-al}), then, for
$\alpha=\{\alpha_k, \, k\geq 1\}$ and $\beta=\{\beta_k,\, k\geq 1\}$ we have
\bel{WProd-xial}
\xi_{\alpha}\diamond \xi_{\beta}=\sqrt{\frac{(\alpha+\beta)!}{\alpha!\beta!}}\,
\xi_{\alpha+\beta},
\ee
where $\alpha+\beta=\{\alpha_k+\beta_k,\, k\geq 1\}$ and
$\alpha!=\prod_{k\geq 1}\alpha_k!=\alpha_1!\alpha_2!\alpha_3!\cdots$.
Using (\ref{WProd-xial}) and linearity, we now define the Wick product of two
arbitrary elements of $\HX$,
\bel{WP-def}
\left(\sum_{\alpha\in \cI} c_{\alpha}\xi_{\alpha}\right) \diamond
\left(\sum_{\beta\in \cI} d_{\beta}\xi_{\beta}\right)=
\sum_{\alpha,\beta\in \cI}c_{\alpha}d_{\beta}
\sqrt{\frac{(\alpha+\beta)!}{\alpha!\beta!}}\, \xi_{\alpha+\beta},
\ee
as long as the series on the right hand side converges in $\HX$.
In general, there is no guarantee that, for  $\xi,\eta\in \HX$, the
Wick product $\xi\diamond \eta $ belongs to $\HX$. For example, let
$\xi\in \HXt$, $\bE \xi^2=1$,  and $\eta=\sum_{n=1}^{\infty}H_n(\xi)/(n\sqrt{n!})$.
Then, treating $\xi$ as the first element of the orthonormal basis
in $\HXt$, we have $\xi=\xi_1=\xi_{\epsilon_1}$ and
$\eta=\sum_{n\geq 1}n^{-1}\xi_{n\epsilon_1}$. Then, by (\ref{WProd-xial}),
$\xi_1\diamond\xi_{n\epsilon_1}=\sqrt{n+1}\xi_{(n+1)\epsilon_1}$, so that
$$
\xi\diamond \eta=\sum_{n=1}^{\infty}
\frac{\sqrt{n+1}}{n}\ \xi_{(n+1)\epsilon_1},
$$
and the series does not converge in $\HX$.

Note that, unlike the usual product, the Wick product of two random variables must be
computed using the chaos expansion (\ref{WP-def}).
The lack of an easy criterion for the convergence in (\ref{WP-def}) is one reason
for considering weighted chaos spaces.
 In the case when $\mfX$ is a white noise
over $L_2(\bR^n)$, weighted chaos spaces are described, for example, in the
books \cite{HKPS} and \cite{HOUZ} (see also \cite{LR1}). The extension of these
spaces to other Gaussian fields is straightforward, but is outside the
scope of our discussion.

Let us summarize the main properties of the Wick product:
\begin{itemize}
\item $\xi\diamond \eta=\eta\diamond \xi$;
\item $\xi\diamond(\eta\diamond \zeta)=(\xi\diamond \eta)\diamond \zeta$;
\item $\xi\diamond(\eta+\zeta)=\xi\diamond \eta+\xi\diamond \zeta$;
\item $\xi\diamond \eta=\xi\,\eta$ if $\xi, \eta\in \HXt$ and $\bE(\xi\eta)=0$.
\item $\xi\diamond \eta=\xi\,\eta$ if either $\xi$ or $\eta$ is an element
of $\HX^{(0)}$, that is, non-random.
\end{itemize}

Similar to ordinary powers, we define Wick powers of a random variable
$\eta\in \HX$: $\eta^{\diamond n}=\eta\diamond\cdots\diamond \eta$.
Replacing ordinary powers with Wick powers in a Taylor series for a function
$f$  leads
to the notion of a Wick function $f^{\diamond}$. For example, the Wick
exponential $e^{\diamond \eta}$ is defined by
\bel{Wick-exp}
e^{\diamond \eta}=\sum_{n=1}^{\infty} \frac{\eta^{\diamond n}}{n!}
\ee
and satisfies $e^{\diamond (\xi+\eta)}=e^{\diamond \xi}\diamond e^{\diamond \eta}.$
If $\eta\in\HXt$, then direct computations show that
\bel{WickExp}
e^{\diamond \eta}=e^{\eta-\frac{1}{2}\bE\eta^2}.
\ee
For more information on the Wick functions, see \cite{HOUZ}.

Just as the chaos decomposition (\ref{OrthDec}), the Wick product doest not
depend on the choice of the orthonormal basis in $\HXt$. For example, if
$\eta_1, \ldots, \eta_k$ are elements of $\HXt$ and
$m_1, \ldots, m_k$ are non-negative integers, then
$
\eta_1^{\diamond m_1}\diamond \cdots \diamond  \eta_k^{\diamond m_k}
$
is the orthogonal projection of $\prod_{j=1}^k\eta_j^{m_j}$ onto
$\HX^{(N)}$, where $N=m_1+\cdots+m_k$. For more details, we refer to
\cite{HKPS,HOUZ,Maj}.

\section{Stochastic Integration}
\label{sec:SI}

In the definition of a generalized random field $\mfX$ over a Hilbert space $\bH$,
 we consider random variables
$\mfX(f)$ for non-random $f\in \bH$. In this section, we define $\mfX(\eta)$
for $\bH$-valued random elements $\eta$.

As a motivation, consider a white noise $\mfB$ over $L_2((0,T))$.
By Theorem \ref{th:wn11},
$W(t)=\mfB(\chi_t)$ is a standard Brownian motion; according to
(\ref{wn:wp}),
\bel{SI-aux1}
W(t)=\sum_{k=1}^{\infty} M_k(t)\xi_k,
\ee
 where $\xi_k=\mfB(m_k)$,
$M_k(t)=\int_0^tm_k(s)ds$, and $\{m_k,\,k\geq 1\}$ is
an orthonormal basis in $L_2((0,T))$.
 Being a continuous function, $W=W(t)$ is an
element of $L_2((0,T))$.   To define $\mfB$ on $W$ using (\ref{SI-aux1}),
  one possibility is
to set $\mfB^{\circ}(W)=\sum_{k=1}^{\infty} \xi_k\mfB(M_k);$
then direct computations show that $\sum_{k=1}^{\infty} \xi_k\mfB(M_k)=W^2(T)/2$.
In other  words, $\mfB^{\circ}(W)=\int_0^TW(t)\circ dW(t)$,
where $\circ$ denotes the Stratonovich integral.
Another possibility is to set
$\mfB^{\diamond}(W)=\sum_{k=1}^{\infty} \xi_k\diamond\mfB(M_k);$
then direct computations show that $\sum_{k=1}^{\infty} \xi_k\diamond\mfB(M_k)
=(\mfB(T)^2-T)/2$.
In other  words, $\mfB^{\diamond}(W)=\int_0^T W(t) dW(t)$,
the \Ito~integral.

We will now use this example to define stochastic integrals
with respect to a white noise $\mfB$  over a
separable Hilbert space $\bH$. Let $\{m_k,\ k\geq 1\}$ be an
orthonormal basis in $\bH$. Define $\xi_k=\mfB(m_k)$ and
 $\xi_{\alpha},\, \alpha \in \cI,$ according to (\ref{xi-al}).

 \begin{definition} An  $\bH$-valued random
element $\eta$ is called \und{$(\mfB,\bH)$-admissible} if
$\bE\|\eta\|_{\bH}^2<\infty$ and,
for every $f\in \bH$, the random variable $(\eta,f)_{\bH}$ is
$\cF^{\mfB}$-measurable.
\end{definition}

 By Theorem \ref{th:basis} and Proposition \ref{prop:wne},
 every $(\mfB,\bH)$-admissible $\eta$ has chaos expansion
\bel{rep-eta}
\eta=\sum_{\alpha\in \cJ}\eta_{\alpha} \xi_{\alpha}, \
\eta_{\alpha}=\bE(\eta\xi_{\alpha})\in \bH.
\ee

\begin{definition}\label{def:main}
Let $\eta$ be  $(\mfB,\bH)$-admissible with chaos expansion (\ref{rep-eta}).\\
The \und{\Ito~ stochastic integral} of $\eta$ with respect to $\mfB$ is
\bel{def:itoi}
\mfB^{\diamond}(\eta)=\sum_{\alpha\in \cI} \mfB(\eta_{\alpha})\diamond \xi_{\alpha},
\ee
where $\diamond$ is the Wick product.
The \und{Stratonovich stochastic integral} of $\eta$ with respect to $\mfB$ is
\bel{def:strati}
\mfB^{\circ}(\eta)=\sum_{\alpha\in \cI} \mfB(\eta_{\alpha})\cdot\xi_{\alpha},
\ee
where $\cdot$ is the usual product.
\end{definition}

Since every generalized Gaussian field and every Gaussian process can be
represented as a white noise over a suitable Hilbert space,  formulas
(\ref{def:itoi}) and (\ref{def:strati}) define stochastic integral
with respect to any Gaussian process or field. We will see below that
these formulas also
provide a chaos expansion of the integral in terms of the chaos expansion
of the integrand; note that neither (\ref{def:itoi}) nor (\ref{def:strati})
is a chaos expansion in the sense of (\ref{rep-eta}).
 The two immediate question that are raised by the above definition
 and will be discussed below are
(a) the convergence of the series, and (b)
the dependence of the integrals on the choice of the basis in $\bH$.

We start by deriving the chaos expansion of the integrals without
investigating the question of convergence.

\begin{theorem}\label{th:st-in-ce}
Let $\eta$ be $(\mfB,\bH)$-admissible with chaos expansion (\ref{rep-eta}),
and assume that
\bel{rep-eta-al}
\eta_{\alpha}=\sum_{k=1}^{\infty} \eta_{\alpha,k}m_k.
\ee
Then
\bel{ito-chaos}
\mfB^{\diamond}(\eta)=\sum_{\alpha\in \cI}\left( \sum_{k=1}^{\infty}
\sqrt{\alpha_k}\eta_{\alpha-\epsilon_k,k}\right) \xi_{\alpha},
\ee
\bel{strat-chaos}
\mfB^{\circ}(\eta)=
\sum_{\alpha\in \cI}
\left( \sum_{k=1}^{\infty}\left(
\sqrt{\alpha_k}\eta_{\alpha-\epsilon_k,k}+
\sqrt{\alpha_k+1}\eta_{\alpha+\epsilon_k,k}\right)\right) \xi_{\alpha}.
\ee
\end{theorem}

\pr By (\ref{rep-eta-al}) and linearity,
$$
\mfB(\eta_{\alpha})=\sum_{k=1}^{\infty}\eta_{\alpha,k}\mfB(m_k)=
\sum_{k=1}^{\infty} \eta_{\alpha,k}\xi_k.
$$
Therefore,
\bel{ito-aux}
\mfB^{\diamond}(\eta)=\sum_{\alpha\in \cI}
\sum_{k=1}^{\infty} \eta_{\alpha,k}\xi_k\diamond \xi_{\alpha}=
\sum_{\alpha\in \cI}\sum_{k=1}^{\infty}
\sqrt{\alpha_k+1}\,\eta_{\alpha,k} \, \xi_{\alpha+\epsilon_k},
\ee
where the last equality follows from (\ref{WProd-xial}); recall
that $\epsilon_k$  is the multi-index with the only non-zero entry, equal to one,
at position $k$. By shifting the summation index, we get (\ref{ito-chaos})
Note that, for every $\alpha \in \cI$, the inner sum in (\ref{ito-chaos})
contains finitely many non-zero terms.

To establish (\ref{strat-chaos}), we write,
similar to (\ref{ito-aux}),
$$
\mfB^{\circ}(\eta)=\sum_{\alpha\in \cI}
\sum_{k=1}^{\infty} \eta_{\alpha,k}\,\xi_k \xi_{\alpha},
$$
 and, instead of (\ref{WProd-xial}), use the following
property of the Hermite polynomials,
$$
H_1(x)H_n(x)=H_{n+1}(x)+nH_{n-1}(x),
$$
which implies
\bel{strat-chaos0}
\xi_k\xi_{\alpha}=\left(\prod_{j\not=k}\frac{H_{\alpha_j}(\xi_j)}{\sqrt{\alpha_j!}}
\right)\frac{H_1(\xi_k)H_{\alpha_k}(\xi_k)}{\sqrt{\alpha_k!}}=
\sqrt{\alpha_k+1}\, \xi_{\alpha+\epsilon_k}+
\sqrt{\alpha_k}\xi_{\alpha-\epsilon_k},
\ee
and then (\ref{strat-chaos}) follows.
\endproof

Now, let us address the questions of convergence and independence of basis.
The Cauchy-Schwartz inequality implies that if
\bel{D12}
\sum_{\alpha\in \cI} |\alpha|\, \|\eta_{\alpha}\|^2_{\bH}<\infty,
\ee
 then
$\mfB^{\diamond}(\eta)\in \HB$. Further examination of (\ref{ito-chaos})
shows that, for every $(\mfB,\bH)$-admissible $\eta$ satisfying (\ref{D12}),
$\mfB^{\diamond}(\eta)$ coincides with the action of the divergence
operator (adjoint of the Malliavin derivative, see \cite{Nualart})
 on $\eta$ and therefore
does not depend on any arbitrary choices, such as the basis in $\bH$.
In particular, if $\bH=L_2(I)$, then $\mfB^{\diamond}(\eta)$ is the
\Ito-Skorokhod integral of $\eta$.
On the other hand, (\ref{ito-chaos}) allows the extension of $\mfB^{\diamond}$
to weighted chaos spaces, similar to those considered in \cite{HKPS,HOUZ,LR1}.

For the Stratonovich integral $\mfB^{\circ}(\eta)$,
 note that the Mallivain derivative $\bD$ of $\xi_{\alpha}$
satisfies
$$
\bD \xi_{\alpha}=\sum_{k=1}^{\infty}\sqrt{\alpha_k}\xi_{\alpha-\epsilon_k}m_k;
$$
this follows directly from the definition of $\bD$ \cite[Definition 1.2.1]{Nualart}
and the relation $H_n'(x)=nH_{n-1}(x)$.
As a result, we use (\ref{strat-chaos0}) to re-write  (\ref{strat-chaos}) as
\bel{strat-chaos1}
\mfB^{\circ}(\eta)=\mfB^{\diamond}(\eta)+\sum_{\alpha\in \cI}
(\eta_{\alpha},\bD \xi_{\alpha})_{\bH}.
\ee
In particular, if $\mfB$ is a white noise over $L_2(I)$ and $\eta=\eta(t)$
is in the domain of the Malliavin derivative, then
\bel{strat-chaos2}
\mfB^{\circ}(\eta)=\mfB^{\diamond}(\eta)+\int_I\bD_t \eta dt,
\ee
where
$$
\bD_t \eta=\sum_{\alpha\in \cI}\eta_{\alpha}(t)
\left(\sum_{k=1}^{\infty}\sqrt{\alpha_k}\,\xi_{\alpha-\epsilon_k}m_k(t)
\right);
$$
u1nlike the \Ito~integral, though, condition
  (\ref{D12}) is not enough to ensure the existence of $\mfB^{\circ}(\eta)$
  as an element of $\HX$. When $\bH$ is the Hilbert space of
  functions on an interval $I$, square integrable
   with respect to a (not necessarily Lebesque)
    measure $\mu$, the sufficient conditions for the Stratonovich integrability
    are discussed in \cite[Chpater 3]{Nualart}. Alternatively,
   $\mfB^{\circ}$ can
 be defined in weighted chaos spaces, but the details of the
 construction have yet to be worked out.

 In what follows, we will concentrate on the \Ito~integral.

 Let $\mfX$ be a zero-mean non-degenerate
 generalized Gaussian field over a separable Hilbert space $\bH$.
 As we mentioned earlier, by the second part of Theorem \ref{th:wn},
  $\mfX$ is a white
 noise over a bigger Hilbert space $\bH_{\cR}$, and then
  $\mfX^{\diamond}(\eta)$ can be defined using (\ref{def:itoi}). If the
 space $\bH_{\cR}$ is difficult to describe,
 one can use representation (\ref{wn11}) from
 the first part of Theorem \ref{th:wn} and consider a
 different formula for the stochastic integral:
 \bel{def:itoi111}
 \mfX^{\diamond}(\eta)=\mfB^{\diamond}(\cK^* \eta)
 \ee
 for every $(\mfB, \bH)$-admissible $\eta$. Similar to the non-random integrands,
 the two definition are equivalent if $\mfX$ is non-degenerate.

 Unlike (\ref{def:itoi}), representation (\ref{def:itoi111}) is not
 intrinsic: the operator $\cK^*$ and the white noise $\mfB$ are not
 uniquely determined by $\mfX$. On the other hand, in many examples,
 such as fractional Brownian motion with the Hurst parameter bigger than $1/2$,
 it is possible to take $\bH=L_2(I)$, and then
 (\ref{def:itoi111}) becomes more convenient than (\ref{def:itoi}).
 To derive the chaos expansion of $\mfX^{\diamond}(\eta)$ using
 (\ref{def:itoi111}), fix  an orthonormal basis $\{m_k,\, k\geq 1\}$ in $\bH$,
 define $\xi_k=\mfB(m_k)$,
 and consider the corresponding orthonormal basis
 $\{\xi_{\alpha}, \alpha\in \cI\}$ in $\HB$
  constructed according to (\ref{xi-al}).
  It follows from (\ref{ito-chaos}) that
 \bel{ito-chaos-alt}
\mfX^{\diamond}(\eta)=\sum_{\alpha\in \cI}\left( \sum_{k=1}^{\infty}
\sqrt{\alpha_k}\,\wt{\eta}_{\alpha-\epsilon_k,k}\right) \xi_{\alpha},
\ee
where
$$
 \wt{\eta}_{k,\alpha}=\bE\big((\cK^*\eta,m_k)_{\bH}\,\xi_{\alpha}\big).
 $$

If $H=L_2(I)$, then (\ref{ito-chaos-alt}) becomes
\bel{st-int-25}
\mfX^{\diamond}_t(\eta)=\sum_{\alpha\in \cI}
\left(\sum_{k\geq 1} \sqrt{\alpha_k}
\left(\int_I \eta_{\alpha-\epsilon_k}(t) (\cK m_k)(t) dt\right)\,
\right)\, \xi_{\alpha},
\ee
where $\eta_{\alpha}(t)=\bE\big(\eta(t) \, \xi_{\alpha}\big)$. In this case,
by analogy with the Brownian motion,  $\int_0^t \eta(s)dX(s)$ can be an
alternative notation for $\mfX^{\diamond}_t(\eta)$,
where $X(t)$ is the associated process of $\mfX$.

We conclude this section with a brief discussion of stochastic differential
equations. To introduce the time evolution, we use
the function $\chi_t$, the characteristic function of the interval
$[0,t]$, and define time-dependent stochastic integrals
\bel{time-evol}
\mfB^{\diamond}_t(\eta):=\mfB^{\diamond}(\eta\chi_t),\ \ \ \
\mfX^{\diamond}_t(\eta):=\mfX^{\diamond}(\eta\chi_t).
\ee
These definitions put an obvious restriction on the Hilbert space
$\bH$, which we call Property I: {\em  $\bH$ is a collection of function or
generalized functions and, for every $\eta\in \bH$ and every fixed $t$,
the (point-wise) product $\eta\chi_t$
is defined and belongs to $\bH$.} There is a more significant
restriction on $\bH$, which we illustrate on the following
equation:
\bel{sde00}
u(t)=1+ \mfB^{\diamond}_t(u),\ 0\leq t \leq T,
\ee
where $\mfB$ is white noise over a Hilbert space $\bH$ with Property I.
Let us assume that the solution belongs to $\HB$ so that
$u(t)=\sum_{\alpha \in \cI} u_{\alpha}(t)\xi_{\alpha}$ and
each $u_{\alpha}$ is an element of $\bH$. By (\ref{time-evol}),
we can re-write (\ref{sde00}) as
\bel{sde01}
u(t)=1+\mfB^{\diamond}(u\chi_t),
\ee
and then (\ref{ito-chaos}) implies
\bel{sde02}
u_{\alpha}(t)=1+\sum_{k=1}^{\infty}(u_{\alpha-\epsilon_k}\chi_t,m_k)_{\bH}.
\ee
Thus, the expression $(u_{\alpha-\epsilon_k}\chi_t,m_k)_{\bH}$, as a
function of $t$, must be an element of $\bH$, and the
Hilbert space $\bH$ must have another special property,
which we call Property II: {\em for every
$f,g\in \bH$, the inner product $(f\chi_t,g)_{\bH}$, as a function of
$t$, is an element of $\bH$}. By the Cauchy-Schwartz inequality,
the space
$L_2(I,\mu)$, with $\mu(I)<\infty$, has both Property I and Property II.
 Representation (\ref{def:itoi111}) then allows us to analyze
stochastic equations for certain generalized Gaussian fields. This
analysis should be a subject of a separate paper, and below we
consider only one particular example.

\begin{theorem} \label{th:sode1}
If $\mfX$ is a zero-mean generalized Gaussian field over $L_2((0,T))$, then
the solution of the equation
\bel{sde2}
u(t)=1+\mfX^{\diamond}_t(u)
\ee
is unique in $L_2((0,T); \HX)$ and is given by
\bel{stoch-exp1}
u(t)=e^{\diamond X(t)},
\ee
where $e^{\diamond}$ is the Wick exponential function (\ref{Wick-exp})
 and $X(t)=\mfX(\chi_t)$ is the associated process of $\mfX$.
\end{theorem}

\pr Let $\mfX(f)=\mfB(\cK^* f)$ be a white noise
representation of $\mfX$ over $L_2((0,T))$.
We start by establishing uniqueness of solution
 in $L_2((0,T);\HB)$, which, because of
the inclusion $\HX\subseteq\HB$, is even stronger.
  By linearity, the difference
$Y(t)$ of two solutions of (\ref{sde2}) satisfies $Y(t)=\mfX^{\diamond}_t(Y).$
If $Y(t)=\sum\limits_{\alpha\in \cI} y_{\alpha}(t)\xi_{\alpha}$, then
(\ref{st-int-25}) implies
\bel{s-syst1}
y_{\alpha}(t)=\sum_{k\geq 1}\sqrt{\alpha_k}
\int_0^ty_{\alpha-\epsilon_k}(s)\wt{m}_k(s)ds,
\ee
where $\wt{m}_k=\cK m_k$. In particular, if $|\alpha|=0$, then
$y_{\alpha}(t)=0$ for all $t$. By induction on $|\alpha|$,
$y_{\alpha}(t)=0$ for all $\alpha \in \cI$: if $y_{\alpha}=0$
for all $\alpha$ with $|\alpha|=n$, then, since
$|\alpha-\epsilon_k|=|\alpha|-1$, equality (\ref{s-syst1})
implies $y_{\alpha}=0$ for all  $\alpha$ with $|\alpha|=n+1$.

To establish (\ref{stoch-exp1}), let
$$
\wt{M}_k(t)=\int_0^t (\cK m_k)(s)ds.
$$
By (\ref{asopr11}),
$$
X(t)=\sum_{k=1}^{\infty} \wt{M}_k(t)\xi_k,
$$
and, because of the independence of $\xi_k$ for different $k$,
$$
e^{\diamond X(t)}=\prod_{k\geq 1} e^{\diamond \wt{M}_k(t)\xi_k}=
\sum_{\alpha\in \cI} \frac{\wt{M}^{\alpha}(t)}{\sqrt{\alpha!}}\, \xi_{\alpha},
$$
where
$$
\wt{M}^{\alpha}(t)=\prod_{k=1}^{\infty} \wt{M}_k^{\alpha_k}(t).
$$
Similar to (\ref{s-syst1}), we conclude that if the solution $u=u(t)$ has the
chaos expansion $u(t)=\sum\limits_{\alpha\in \cI}
u_{\alpha}(t)\xi_{\alpha}$, then $u_{\alpha}(t)=1$ if
$|\alpha|=0$ and
\bel{s-syst12}
u_{\alpha}(t)=\sum_{k\geq 1}\sqrt{\alpha_k}
\int_0^tu_{\alpha-\epsilon_k}(s)\wt{m}_k(s)ds,
\ee
if $|\alpha|>0$. Then direct computations show that
$$
u_{\alpha}(t)=\frac{\wt{M}^{\alpha}(t)}{\sqrt{\alpha!}}, \ |\alpha|\geq 1,
$$
satisfies (\ref{s-syst12}):
\begin{equation*}
\begin{split}
\frac{d u_{\alpha}(t)}{dt}&=\frac{1}{\sqrt{\alpha!}}\,\frac{d}{dt}\prod_{k=1}^{\infty}
\wt{M}_k^{\alpha_k}(t)=\frac{1}{\sqrt{\alpha!}}
\sum_{k=1}^{\infty} \alpha_k \wt{M}_k^{\alpha_k-1}(t)\wt{m}_k(t)\prod_{j\not=k}
\wt{M}_j^{\alpha_j}(t)\\
&=\sum_{k=1}^{\infty}\sqrt{\alpha_k}\, \wt{m}_k(t)
 \frac{\wt{M}^{\alpha-\varepsilon_k}(t)}{\sqrt{(\alpha-\varepsilon_k)!}}=
 \sum_{k=1}^{\infty}\sqrt{\alpha_k}\, m_k(t) u_{\alpha-\varepsilon_k}(t).
 \end{split}
 \end{equation*}
\endproof

Theorem \ref{th:sode1} is a generalization of the
familia result that the geometric Brownian motion
$u(t)=e^{W(t)-(t/2)}=e^{\diamond W(t)}$ satisfies
$u(t)=1+\int_0^tu(s)dW(s)$: by  (\ref{stoch-exp1}) and (\ref{WickExp}),
  for a class of zero-mean Gaussian processes $X=X(t)$ with
covariance function $R(t,s)$, and with a suitable interpretation of the
stochastic integral, the solution of the equation
$u(t)=1+\int_0^tu(s)dX(s)$ is
$$
u(t)=e^{X(t)-\frac{1}{2}R(t,t)}.
$$

 The proof of the theorem suggests that stochastic equations
in the  \Ito-Skorokhod sense
are more suitable for analysis using chaos expansion than the equations
in the Stratonovich sense. Indeed, equation (\ref{sde2})
leads to the system of equations (\ref{s-syst12})
 that is   solvable by induction on $|\alpha|$.
By contrast,  equation
$u(t)=1+\mfX^{\circ}_t(u)$ leads to a system that is not
solvable by induction on $|\alpha|$:  according to
 (\ref{strat-chaos}), $u_{\alpha}$ will
depend on both $u_{\alpha-\epsilon_k}$ and $u_{\alpha+\epsilon_k}$.

The arguments used in the proof of Theorem \ref{th:sode1}
 can be extended to more general linear equations
and to generalized fields over $L_2((0,T),\mu)$ for different measures $\mu$,
although the precise results will essentially depend on certain fine properties
of $\mu$.


\begin{thebibliography}{10}

\bibitem{AMN1}
E.~Al\`{o}s, O.~Mazet, and D.~Nualart.
\newblock {Stochastic Calculus With Respect to Fractional Brownian Motion with
  Hurst Parameter Less Than $\frac{1}{2}$}.
\newblock {\em Stochastic Process. Appl.}, 86(1):121--139, 2000.

\bibitem{AMN}
E.~Al\`{o}s, O.~Mazet, and D.~Nualart.
\newblock {Stochastic Calculus With Respect to Gaussian Processes}.
\newblock {\em Ann. Probab.}, 29(2):766--801, 2001.

\bibitem{A}
A.~Amirdjanova.
\newblock {Nonlinear Filtering with Fractional Brownian Motion}.
\newblock {\em Appl. Math. Optim.}, 46(2--3):81--88, 2002.

\bibitem{CM}
R.~H. Cameron and W.~T. Martin.
\newblock {The Orthogonal Development of Nonlinear Functionals in Series of
  {Fourier-Hermite} Functionals}.
\newblock {\em Ann. of Math.}, 48(2):385--392, 1947.

\bibitem{DH}
W.~Dai and C.~C. Heyde.
\newblock {It\^{o}'s Formula with Respect to Fractional Brownian Motion and its
  Application}.
\newblock {\em J. Appl. Math. Stochastic Anal.}, 9(4):439--448, 1996.

\bibitem{DU}
L.~Decreusefond and A.~S. \"{U}st\"{u}nel.
\newblock {Stochasic Analysis of the Fractional Brownian Motion}.
\newblock {\em Potential Anal.}, 10(2):177--214, 1998.

\bibitem{DHP}
T.~E. Duncan, Y.~Hu, and B.~Pasik-Duncan.
\newblock {Stochastic Calculus for Fractional Brownian Motion I: Theory}.
\newblock {\em SIAM J. Control Optim.}, 38(2):582--612, 2000.

\bibitem{DSch2}
N.~Dunford and J.~T. Schwartz.
\newblock {\em {Linear Operators, Part II: Spectral Theory}}.
\newblock Wiley Classics Library Publication, 1988.

\bibitem{HKPS}
T.~Hida, H-H. Kuo, J.~Potthoff, and L.~Sreit.
\newblock {\em {White Noise}}.
\newblock Kluwer, 1993.

\bibitem{HOUZ}
H.~Holden, B.~{\O}ksendal, J.~Ub{\o}e, and T.~Zhang.
\newblock {\em {Stochastic Partial Differential Equations: A Modeling, White
  Noise Functional Approach}}.
\newblock Birkh\"{a}user, 1996.

\bibitem{KBR}
M.~L. Kleptsyna, A.~{Le} Breton, and M.-C. Roubaud.
\newblock {General Approach to Filtering with Fractional Brownian Noises:
  Application to Linear Systems}.
\newblock {\em Stochastics Stochastics Rep.}, 71(1--2):119--140, 2000.

\bibitem{Lin}
S.~J. Lin.
\newblock {Stochastic Analysis of Fractional Brownian Motions}.
\newblock {\em Stochastics Stochastics Rep.}, 55(1--2):121--140, 1995.

\bibitem{LR1}
S.~V. Lototsky and B.~L. Rozovskii.
\newblock {Wiener Chaos Solutions of Linear Stochastic Evolution Equations}.
\newblock {\em Ann. Probab.}, 34(2):638--662, 2006.

\bibitem{Maj}
P.~Major.
\newblock {\em {Multiple Wiener-It{\^{o}} Integrals. With Applications to Limit
  Theorems}}, volume 849 of {\em {Lecture Notes in Mathematics}}.
\newblock Springer, 1981.

\bibitem{Nualart}
D.~Nualart.
\newblock {\em {The Malliavin Calculus and Related Topics, 2nd Edition}}.
\newblock Springer, 2006.

\bibitem{PT}
V.~Pipiras and M.~S. Taqqu.
\newblock {Integration Questions Related to Fractional Brownian Motion}.
\newblock {\em Probab. Theory Related Fields}, 118(2):251--291, 2000.

\end{thebibliography}

\end{document}